\documentclass{imsart}

\usepackage{amsmath}
\usepackage{bm}
\usepackage[numbers]{natbib}
\usepackage[colorlinks,citecolor=blue,urlcolor=blue]{hyperref}
\usepackage{cleveref,booktabs}
\usepackage{graphicx}

\startlocaldefs
\endlocaldefs

\begin{document}

\begin{frontmatter}

\title{On testing mean proportionality of multivariate normal
  variables}
\runtitle{Testing mean proportionality}


\author{\fnms{Etaash} \snm{Katiyar}\ead[label=e1]{ek493@cam.ac.uk}}
\and
\author{\fnms{Qingyuan} \snm{Zhao}\corref{}\ead[label=e2]{qyzhao@statslab.cam.ac.uk}}

\address{Centre for Mathematical Sciences\\
  Wilberforce Road, Cambridge CB3 0WB \\
  United Kingdom \\
\printead{e1,e2}}

\runauthor{E. Katiyar \& Q. Zhao}

\affiliation{University of Cambridge}

\begin{abstract}
  This short note considers the problem of testing the null hypothesis that the mean
  values of two multivariate normal variables are proportional. We
  show that the usual likelihood ratio $\chi^2$-test is valid
  non-asymptotically. Our proof relies on expressing the test
  statistic as the minimum eigenvalue of a Wishart variable and using
  a representation of its distribution using Legendre polynomials.
\end{abstract}





\end{frontmatter}


\section{Introduction}
\label{sec:introduction}

Suppose $X$ and $Y$ are independent $p$-dimensional random vectors, $X
\sim \mathrm{N}(\mu_1, I_p)$ and $Y \sim \mathrm{N}(\mu_2,
I_p)$. This paper is concerned with
testing the hypothesis that their mean values are proportional,
$H_0:\mu_1 \propto \mu_2$. That is, we are interested in
testing the hypothesis that there exists a scalar $\eta$ such that
$\mu_2 = \beta \mu_1$. This problem arises naturally in many
applications, such as instrumental variables regression
\citep{zhao2020statistical} and genetic colocalization analysis
\citep{wallace2013statistical}. In fact, our investigation is motivated
by the last application. In genetic colocalization analysis, the
measurements $X$ and $Y$ are the regression coefficients of two phenotypes
on the same genotypes in a genomic region. In practice, they are
usually obtained from different genome-wide association studies. Under
the assumption that the two phenotypes share a single causal genetic
variant in this region, the two sets of regression coefficients should
have proportional means \citep{plagnol2008statistical}.

The mean proportionality testing problem is closely related to
Fieller's theorem and errors-in-variables
regression. Assuming that $H_0$ is indeed true,
\citet{fieller1954some} proposes to construct a confidence interval
for $\beta$ by using the pivot
\begin{equation}
  \label{eq:r-beta}
  R(\beta) = (Y - \beta X)^T (Y - \beta X) / (1 + \beta^2) \sim
  \chi_p^2.
\end{equation}
Because the only unknown quantity in \eqref{eq:r-beta} is $\beta$, one
can obtain an exact confidence interval for $\beta$ by using suitable
quantiles of $\chi_p^2$. However, our interest lies in testing the
existence of such $\beta$ instead of estimating $\beta$ when it is
assumed to exist. Thus, the problem being considered here is a special
case of goodness-of-fit testing for errors-in-variables regression.

Although we have assumed that both $X$ and $Y$ have identity
covariance, the same formulation applies to the more general setting
where $X \sim \mathrm{N}(\mu_1, \Sigma)$ and $Y \sim \mathrm{N}(\mu_2, \Sigma)$,
where the $p \times p$ matrix $V$ is known. In such case, we can
simply consider the transformed variables $\Sigma^{-1/2} X$ and $\Sigma^{-1/2}
Y$, whose mean values, $\Sigma^{-1/2} \mu_1$ and $\Sigma^{-1/2} \mu_2$, are still
proportional under $H_0$. The invariance of $H_0$ under
scaling also means that we can allow the covariance matrix of
$Y$ to be $c\Sigma$, where $c$ can be an unknown scalar.

To test the existence of $\beta$, \citet{plagnol2008statistical}
propose to compare the minimum value of $R(\beta)$ over $\beta$ with
quantiles of $\chi_{p-1}^2$. More specifically, let $\hat{\beta} =
\arg\min_{\beta} R(\beta)$. \citet{plagnol2008statistical} propose to
reject $H_0$ at level $(1-\alpha)$ if $R(\hat{\beta}) >
\chi_{p-1}^2(1-\alpha)$, where $\chi_{p-1}^2(1-\alpha)$ is the
$(1-\alpha)$ quantile of $\chi_{p-1}^2$ for $0 < \alpha < 1$. This test is also
described in \citet{wallace2013statistical} and implemented in a
popular \textsf{R} package called \textsf{coloc}.

To justify the aforementioned $\chi^2$-test,
\citet{plagnol2008statistical} and \citet{wallace2013statistical}
cite asymptotic theory but do not provide a formal argument. It is
straightforward to verify that $R(\beta)$ is twice the negative profile
log-likelihood of $\beta$ under $H_0$ (up to an additive
constant). Moreover, $R(\hat{\beta})$ is exactly twice the
negative log likelihood ratio statistic for testing $H_0:\mu_1
\propto \mu_2$ versus $H_1: \mu_1$ and $\mu_2$ are unrestricted. Because
$R(\beta) \sim \chi^2_p$ and one degree of freedom is
spent on estimating $\beta$, intuitively one may expect that
$R(\hat{\beta})$ converges in distribution to $\chi_{p-1}^2$. However,
this does not immediately follow from Wilk's theorem or the classical
asymptotic theory for likelihood ratio tests, because the dimension of
the parameter space is changing. The null model can be parameterized
by the $(p+1)$ dimensional vector $(\mu_1^T, \beta)$, while the full
model is parameterized by the $2p$ dimensional vector
$(\mu_1^T,\mu_2^T)$. Thus both the null and full model spaces have
growing dimensions. In fact, the very statement $R(\hat{\beta})
\to \chi_{p-1}^2$ in distribution as $p \to \infty$ is not rigorous,
because the distributional limit is changing with $p$. Another
potential concern is that the standard likelihood theory may not apply
if $\mu_1^T\mu_1$ does not grow as fast as $p$ when $p \to \infty$
\citep{zhao2020statistical}.

Notice that $R(\hat{\beta})$ is stochastically dominated by $\chi^2_p$,
meaning that its $(1-\alpha)$ quantile is smaller than the
$\chi^2_p(1-\alpha)$ for all $0 < \alpha < 1$. This result is trivial
because $R(\hat{\beta})
\leq R(\beta)$ by definition and $R(\beta) \sim \chi^2_p$. In the rest
of this article, we will show that the distribution of
$R(\hat{\beta})$ is also stochastically dominated by
$\chi^2_{p-1}$ for all $p \ge 2$.

\section{Exact distribution of $R(\hat{\beta})$}

Our result relies on classical distributional results
on the eigenvalues of a Wishart random variable and is non-asymptotic
(does not require $p \to \infty$). Let $c(\beta) = (\beta, - 1)^T /
\sqrt{1 + \beta^2}$; notice that $c(\beta)$ has $\ell_2$-norm equal to
$1$. Let $S = (X \, Y)^T (X \, Y)$. Observe that
\begin{align*}
  R(\hat{\beta}) &= \inf_{\beta} R(\beta) \\
                 &= \inf_{\beta} (Y - \beta X)^T (Y - \beta X) / (1 + \beta^2) \\
                 &= \inf_{\beta} c(\beta)^T (X \, Y)^T (X \, Y) c(\beta) \\
                 &= \lambda_2(S),
\end{align*}
where $\lambda_1(S) \geq \lambda_2(S) \ge 0$ are the two eigenvalues
of $S$. Because $S$ follows a Wishart distribution (non-central if
$\mu_1 \neq 0$), this allows us to use classical distributional
results on the eigenvalues of a Wishart random variable
\citep{james1964distributions}.

The eigenvalue distribution is much simpler when the Wishart
distribution is central. When the scale matrix is identity,
\citet[Corollary 3.2.19]{muirhead1982aspects} has derived the joint density
function of the eigenvalues. When $S$ is $2 \times 2$, this is given
by
\[
  f_c(\lambda_1,\lambda_2) \propto e^{-(\lambda_1 + \lambda_2)/2}
  (\lambda_1 \lambda_2)^{(p-3)/2} (\lambda_1 - \lambda_2),~\lambda_1 >
  \lambda_2 > 0.
\]
The normalizing constant can be found in \citet{muirhead1982aspects}.

The non-central case is more complicated. \citet[equation
68]{james1964distributions} has given the joint density function of the
eigenvalues of a non-central Wishart random variable. For the case we
are considering ($S$ is $2 \times 2$ and scale matrix is identity),
the joint density is given by, for $\lambda_1 > \lambda_2
> 0$,
\begin{equation}
  \label{eq:joint-density}
  \begin{split}
    &f(\lambda_1,\lambda_2) \propto \\
    & {}_0 F_1\left(\frac{p}{2};
    \begin{pmatrix}
      (1+\beta^2) \mu_1^T \mu_1 / 4 & 0 \\
      0 & 0 \\
    \end{pmatrix},
    \begin{pmatrix}
      \lambda_1 & 0 \\
      0 & \lambda_2
    \end{pmatrix}
  \right) e^{- (\lambda_1 + \lambda_2)/2} (\lambda_1
  \lambda_2)^{(p-3)/2} (\lambda_1 - \lambda_2),
  \end{split}
\end{equation}
where ${}_0 F_1$ is the generalized hypergeometric function of two
matrix arguments defined in, for instance,
\citet{james1964distributions}. The normalizing constant for
$f(\lambda_1,\lambda_2)$ can also be found
there. By definition, $f(\lambda_1,\lambda_2)$ reduces to
$f_c(\lambda_1,\lambda_2)$ when $\mu_1 = 0$. By using equation (1.13)
in \citet{muirhead1975expressions}, we
can write ${}_0 F_1$ as a series
\begin{equation}
  \label{eq:0F1}
  \begin{split}
  &{}_0 F_1\left(\frac{p}{2};
    \begin{pmatrix}
      (1+\beta^2) \mu_1^T \mu_1 / 4 & 0 \\
      0 & 0 \\
    \end{pmatrix}, \begin{pmatrix}
      \lambda_1 & 0 \\
      0 & \lambda_2
    \end{pmatrix}
  \right) \\
  =& \sum_{j=0}^{\infty}
  \frac{\left(\left(1+\beta^2\right)\mu_1^T\mu_1\right)^j\left(\lambda_1
      \lambda_2\right)^{j/2}}{4^j (p/2)_j
    j!} P_j\left(\frac{\lambda_1 + \lambda_2}{2 \sqrt{\lambda_1
        \lambda_2}}\right),
  \end{split}
\end{equation}
where $P_j$ is the Legendre polynomial of degree $j$ and $(p/2)_j$ is
the rising factorial $(p/2)_j = (p/2) (p/2+1) \cdots (p/2 + j-1)$. The
Legendre polynomials can be found using the initial polynomials
$P_0(x) = 0$, $P_1(x) = x$, and Bonnet's recursion formula
\begin{equation}
  \label{eq:bonnet-recursion}
  (j+1) P_{j+1}(x) = (2j+1) x P_j(x) - j P_{j-1}(x), \quad j=1,2,\dotsc.
\end{equation}
The derivative of the Legendre polynomials can be computed by
\begin{equation}
  \label{eq:bonnet-recursion-2}
  \frac{\text{d}}{\text{d} x} P_j(x) = \frac{j xP_{j}(x) - j
    P_{j-1}(x)}{x^2-1}, \quad j=1,2,\dotsc.
\end{equation}
These formulae will be quite useful in our proof below.

There has been a lot of investigations on the distribution of the
largest eigenvalue of a sample covariance matrix, as it is closely
related to selecting the number of principal components in a principal
component analysis \citep[see
e.g.,][]{johnstone2001distribution}. There has been relatively less
interest in studying the distribution of the smallest
eigenvalue. \citet{edelman1991distribution} has derived a recursion for
the distribution of the smallest eigenvalue of a central Wishart
variable without resorting to zonal polynomials and hypergeometric
functions of matrix arguments. Although that formula can be used to
efficiently compute the distribution numerically, we could not use it
to prove that $\lambda_2(S)$ is stochastically dominated by
$\chi^2_{p-1}$. Instead, in the proof below we will directly obtain
the distribution of $\lambda_2(S)$ by integrating
$f(\lambda_1,\lambda_2)$ over $\lambda_1$.

\section{Proof of stochastic dominance}

Consider two cumulative distribution functions $F(x)$ and
$G(x)$. Suppose they are defined on the same support (in our case,
$(0, \infty)$) and have density functions $f(x)$ and $g(x)$,
respectively. The distribution $F$ is said to have (first-order)
stochastic dominance over $G$ if $F(x) \leq G(x)$ for all $x$. To show
stochastic dominance, it is sufficient to establish monotone
likelihood ratio property, that is, $f(x)/g(x)$ is increasing in $x$
\citep[Proposition 4.3]{wolfstetter1999topics}.

We first consider the central case of our problem, as the distribution of
$\lambda_2(S)$ is much simpler. In this case, the marginal density
function of $\lambda_2(S)$ is given by
\begin{align*}
  f_c(\lambda_2)& \propto \int_{\lambda_2}^{\infty}e^{-(\lambda_1+\lambda_2)/2} \left(\lambda_1 \lambda_2\right)^{(p-3)/2}\left(\lambda_1-\lambda_2\right)\text{d}\lambda_1 \\
                & \propto e^{-\lambda_2/2}\lambda_2^{(p-3)/2}\left[\Gamma\left(\frac{p+1}{2},\frac{\lambda_2}{2}\right) - \frac{\lambda_2}{2} \Gamma\left(\frac{p-1}{2},\frac{\lambda_2}{2}\right) \right],
\end{align*}
where $\Gamma$ is the incomplete gamma function $\Gamma(s,x)=
\int_x^\infty t^{s-1}e^{-t} \text{d}t$. To show that $\lambda_2(S)$ is
stochastically dominated by $\chi^2_{p-1}$, whose density function is
proportional to $e^{-\lambda_2/2} \lambda_2^{(p-3)/2}$, it suffices to
show that the likelihood ratio (up to a multiplicative constant)
\[
  g_c(\lambda_2) = \Gamma\left(\frac{p+1}{2},\frac{\lambda_2}{2}\right) -
  \frac{\lambda_2}{2} \Gamma\left(\frac{p-1}{2},\frac{\lambda_2}{2}\right)
\]
is decreasing in $\lambda_2$. This can be verified by checking that
the derivative of $g_c(\lambda_2)$ is equal to $- \Gamma((p-1)/2,
\lambda_2/2) < 0$, using $(\partial / \partial x)
\Gamma(s,x) = - x^{s-1} e^{-x}$.

For the non-central case, the marginal density of $\lambda_2$
\[
  f(\lambda_2) = \int_{\lambda_2}^{\infty} f(\lambda_1,\lambda_2)
  \text{d} \lambda_1
\]
is much more complicated. A natural idea is to show that the
distribution in the non-central case is dominated by the central
case. However, this is not true (see next section).

Motivated by the proof for the central case, it suffices to show that
the likelihood ratio of $\lambda_2(S)$ to \(\chi_{q-1}^{2}\) is
decreasing, that is,
\[
  \frac{\partial}{\partial\lambda_2}
  \left(\frac{f(\lambda_2)}{e^{-\lambda_2/2}\lambda_2^{(p-3)/2}}
  \right) \leq 0.
\]
By using \eqref{eq:joint-density}, \eqref{eq:0F1}, and Leibniz's rule
for differentiation, this is equivalent to showing
\[
  \begin{split}
  \int_{\lambda_2}^{\infty} &e^{-\lambda_1/2} \lambda_1^{(p-3)/2} \cdot \\
  &  \cdot \frac{\partial}{\partial \lambda_2}
  \sum_{j=0}^{\infty} (\lambda_1 - \lambda_2)
  \frac{\left(\left(1+\beta^2\right)\mu_1^T\mu_1\right)^j\left(\lambda_1
      \lambda_2\right)^{j/2}}{4^j (p/2)_j
    j!} P_j\left(\frac{\lambda_1 + \lambda_2}{2 \sqrt{\lambda_1
        \lambda_2}}\right) \text{d} \lambda_1\leq 0.
  \end{split}
\]
Thus, it suffices to show that for all $\lambda_1 > \lambda_2$,
\[
  \frac{\partial}{\partial \lambda_2}  (\lambda_1 - \lambda_2)
  \lambda_2^{j/2}
  P_j\left(\frac{\lambda_1 + \lambda_2}{2 \sqrt{\lambda_1
        \lambda_2}}\right) \leq 0, \quad j = 0,1,\dotsc.
\]
This trivially holds for $j = 0$. Let $x = (\lambda_1 + \lambda_2)/(2
\sqrt{\lambda_1 \lambda_2})$, so it suffices to show the following
inequality for $j \geq 1$:
\[
  \left\{ (\lambda_1/\lambda_2) (j/2) - (j/2 + 1) \right \} P_j(x) +
  (\lambda_1 - \lambda_2) \frac{\text{d}}{\text{d} x}P_j(x)
  \frac{\partial x}{\partial \lambda_2} \leq 0.
\]
By using $\partial x / \partial \lambda_2 = (\lambda_2 - \lambda_1) /
(4 \sqrt{\lambda_1 \lambda_2^3})$, this is equivalent to
\[
  \left\{\lambda_1 (j/2) - \lambda_2 (j/2 + 1) \right \} P_j(x) \leq
  \frac{(\lambda_1 - \lambda_2)^2}{4 \sqrt{\lambda_1 \lambda_2}}
  \frac{\text{d}}{\text{d} x}P_j(x).
\]
By using \eqref{eq:bonnet-recursion-2} and $x - \sqrt{x^2 - 1} =
\sqrt{\lambda_2/\lambda_1}$, this can be shown to be equivalent to
\begin{equation}
  \label{eq:key-inequality}
  \frac{P_{j-1}(x)}{P_j(x)} \leq
  \frac{j+1}{j} \left(  x - \sqrt{x^2 - 1} \right)\quad \text{for all
  } x > 1.
\end{equation}
In summary, we have reduced the proof of stochastic dominance to
showing the above inequality for Legendre polynomials.

We prove the inequality \eqref{eq:key-inequality} using induction. It
is easy to check that \eqref{eq:key-inequality} holds for $j=1$. Next
we assume \eqref{eq:key-inequality} holds for $j$. By Bonnet's
recursion \eqref{eq:bonnet-recursion} and plugging in
\eqref{eq:key-inequality}, we get
\begin{align*}
  \frac{P_{j+1}(x)}{P_{j}(x)} &= \frac{2j+1}{j+1}x -
                                \frac{j}{j+1}\frac{P_{j-1}(x)}{P_j(x)} \\
                              &\geq \frac{2j+1}{j+1}x - \left(  x - \sqrt{x^2 - 1} \right) \\
                              &= \frac{j}{j+1}\left(x + \sqrt{x^2 - 1}\right) + \frac{1}{j+1}
                                \sqrt{x^2 - 1} \\
\end{align*}
By using $(x + \sqrt{x^2 - 1})(x - \sqrt{x^2 - 1}) = 1$, this can be
written as
\[
  \frac{P_{j+1}(x)}{P_{j}(x)} - \frac{j+1}{j+2}\frac{1}{x - \sqrt{x^2
      - 1}} \geq 
  \frac{(j+1) \sqrt{x^2 - 1} -
    x}{(j+1)(j+2)}. \\
\]
The right hand side is non-negative if and only if $x \geq (j+1) /
\sqrt{j (j+2)}$. Therefore, we cannot rely on using the induction
hypothesis to prove
\begin{equation} \label{eq:inequality-induction}
  \frac{P_j(x)}{P_{j+1}(x)} \leq \frac{j+2}{j+1} \left(  x - \sqrt{x^2
      - 1} \right)
\end{equation}
for $1 < x < (j+1) / \sqrt{j (j+2)}$. Instead, we prove this
inequality directly for such $x$ by noticing that $P_{j-1}(x) \leq
P_j(x)$ for all $x \geq 1$ and $j = 1,2,\dotsc$ (this can be shown by
induction). Thus
\begin{align*}
  \frac{P_{j+1}(x)}{P_j(x)} =  \frac{2j+1}{j+1}x -
  \frac{j}{j+1}\frac{P_{j-1}(x)}{P_j(x)}
  \geq \frac{2j+1}{j+1}x - \frac{j}{j+1}
  \geq \frac{2j+1}{j+1}x - \frac{j}{j+1}x  = x,
\end{align*}
Notice that for $x < (j+1) / \sqrt{j (j+2)}$, we have
\begin{align*}
  \frac{\left(\frac{j+1}{j+2} \right) \left(x + \sqrt{x^2-1} \,\,
  \right) }{x} & =
                 \left(\frac{j+1}{j+2}\right)\left(1+\sqrt{1-\frac{1}{x^2}}\right)
                 < 1.
\end{align*}
It is then straightforward to verify \eqref{eq:inequality-induction}.

\section{Numerical illustration}
\label{sec:numer-illustr}

\begin{figure}[htbp]
  \centering
  \includegraphics[width = \textwidth]{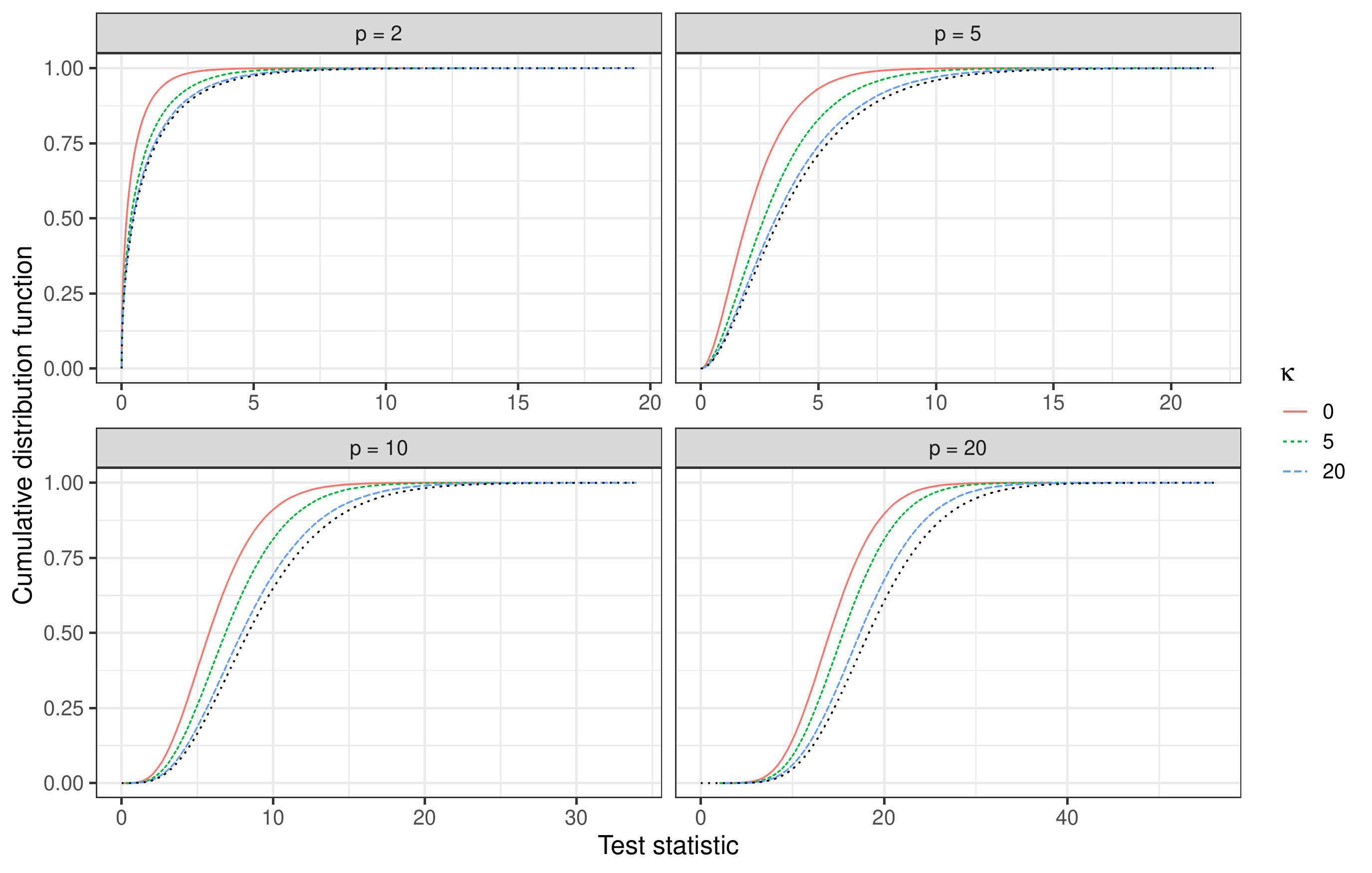}
  \caption{Simulated distribution function of the test statistic
    $R(\hat{\beta})$. The different curves correspond to different
    values of the noncentrality parameter $\kappa = (1+\beta^2) \mu_1^T \mu_1$;
    the rightmost dotted curve in each panel is the distribution
    function of $\chi^2_{p-1}$.}
  \label{fig:cdf}
\end{figure}

\Cref{fig:cdf} shows the simulated distribution distribution of the
test statistic $R(\hat{\beta})$ in various settings. Apart from the
dotted curve which corresponds to $\chi^2_{p-1}$, each curve is
obtained from 100,000 simulations. Notice that the
distribution only depends on $\beta$ and $\mu_1$ through $\kappa = (1
+ \beta^2) \mu_1^T \mu_1$, which can be shown using equation 68 in
\citet{james1964distributions}. We will call $\kappa$ the
noncentrality parameter. From this figure, it appears that the
distribution function of $R(\hat{\beta})$ becomes smaller and
approaches $\chi^2_{p-1}$ as $\kappa$ increases. However, for a fixed
$\kappa$, $\chi^2_{p-1}$ becomes a worse approximation as $p$
increases. \Cref{tab:size} lists the simulated size (type I error) of the
$\chi^2$-test in various settings. The test becomes more conservative
as the noncentrality parameter $\kappa$ decreases and the dimension
$p$ increases.

\begin{table}[htbp]
  \centering
  \caption{Simulated size (type I error) of the $\chi^2$-test for
    different dimension $p$,
    significance level $\alpha$, and noncentrality parameter $\kappa =
    (1 +
    \beta^2) \mu_1^T \mu_1$.}
  \label{tab:size}
  \begin{tabular}{l|ccc|ccc|ccc}
    \toprule
    & \multicolumn{3}{c|}{$\alpha = 1\%$} & \multicolumn{3}{c|}{$\alpha =
                                             5\%$} &
                                                      \multicolumn{3}{c}{$\alpha
                                                      = 10\%$} \\
    & $\kappa=0$ & $\kappa=5$ & $\kappa=20$ & $\kappa=0$ & $\kappa=5$ &
                                                                        $\kappa=20$ & $\kappa=0$ & $\kappa=5$ & $\kappa=20$ \\
    \midrule
    $p=2$  & 0\% & 0.2\% & 0.7\% & 0.3\% & 2.2\% & 4.3\% & 1.3\% & 5.8\% & 9\% \\
    $p=5$  & 0\% & 0.1\% & 0.6\% & 0.2\% & 1.3\% & 3.8\% & 0.8\% & 3.7\% & 8.1\% \\
    $p=10$  & 0\% & 0.1\% & 0.4\% & 0.1\% & 0.8\% & 3.2\% & 0.7\% & 2.5\% & 7.2\% \\
    $p=20$  & 0\% & 0\% & 0.3\% & 0.1\% & 0.5\% & 2.5\% & 0.5\% & 1.7\% & 6\% \\
    \bottomrule
  \end{tabular}
\end{table}

\bibliographystyle{imsart-nameyear}
\bibliography{ref}




\end{document}